\documentclass[a4paper]{article}
 \pdfoutput=1          
 \UseRawInputEncoding  
\usepackage[textwidth=14.5cm]{geometry}  
\usepackage{blindtext}                   
\usepackage{amsmath,amsfonts,amssymb,amsthm}
\usepackage{latexsym}
\usepackage{enumerate}
\usepackage{graphics}
\usepackage{epsfig}
\usepackage{mathdots}
\usepackage{caption}
\usepackage{subcaption}
\usepackage{setspace}
\usepackage{multirow}
\usepackage{dcolumn}
\usepackage{url}
\usepackage{tikz,pgfplots}
\usetikzlibrary{arrows,backgrounds}
\usepackage{ragged2e}
\usepackage[numbers,sort&compress]{natbib}

\newtheorem{thm}{Theorem}
\newtheorem{lem}{Lemma}

\newtheorem{defn}{Definition}  
\newtheorem{prop}{Proposition}  
\newtheorem{exmp}{Example}  
\newenvironment{solution}{\emph{Solution.}}  

\author{Yongwen Zhu\\
\small School of Mathematics and Information Science, Yantai University,\\Yantai City 264005, P.R. China
\\ \small Email: {zyw@ytu.edu.cn}
}

\title{The Plum-Blossom Wedge Product Method for Multiplication and Division of Multi-digit Integers
}

\date{}

\begin{document}

\begin{sloppypar}     

\maketitle

\abstract{
For two integers $a$ and $b$, the plum-blossom product of $a$ and $b$ is defined as the ones of the usual product of $a$ and $b$ if the ones is less than or equal to $3$, or otherwise as the ones minus $10$.
Under this operation, all integers is a commutative semigroup. Since some fine properties of this semigroup, the plum-blossom product method is the efficient method of multiplication and division of long integers, can be used to mental arithmetic and to computer science.
In this paper, we improve the theory of plum-blossom products  so that the new concept--the  plum-blossom wedge product is proposed. The plum-blossom wedge product is also of some good features so that it can be applied to multiplication and division of long integers.
The plum-blossom wedge product method of multiplication and division is not only efficient in mental arithmetic, but also could applied to improve the design of multipliers and dividers of computer.
}



\section{Introduction}\label{introduction}

As we have known, multipliers requiring large bit lengths have a major impact on the performance of many applications, such as cryptography, digital signal processing (DSP) and image processing. As multipliers take a long time for execution so there is a need of fast multiplier to save the execution time. Vedic algorithms are used to the design of computer processors for enhancing speed and performance, see
\cite{Bansal,Garg,Gupta,Gurumurthy,Kavita,Kayal,Mathur,Paramasivam,Pradhan,Ramalatha,Rashno,Sahu}.
For more literatures involving in the hardware architecture, the fully homomorphic encryption, the combined single trace attack on global shuffling long integer multiplication and the configurable long integer multiplier, etc., one may see \cite{Rafferty,San,Hua,Lee,Shi-Wang}.

In \cite{Harvey-Hoeven}, Harvey and Hoeven presented an algorithm that computes the product of two $n$-bit integers in $O(n \log n)$ bit operations, thus confirming a conjecture of Sch{\"o}nhage and Strassen from 1971.  This new algorithm adopted a novel “Gaussian resampling” technique that reduces the integer multiplication problem to a collection of multidimensional discrete Fourier transforms over the complex numbers, which may be evaluated rapidly by means of Nussbaumer' s fast polynomial transforms.

In his  paper \cite{Zhu2}, the author presented a new multiplication formula, which may be used to design the multiplier but also to mental multiplication.
This algorithm was a generalization of that of Karatsuba and Offman \cite{Rafferty,Harvey-Hoeven,Kavita,Karatsuba}, but it is superior to the latter because the latter is not suitable for oral calculation of multiplication.

Recently, in a series of his research work, see for example \cite{Zhu,Zhu1,Zhu2,Zhu3,Zhu4}, the author proposes some novel method to calculate multiplication and division of long integers, where some new concept such as scissor products and the plum-blossom products are introduced. These method is superior to that of \cite{Shi}. For any two integers, the  plum-blossom product is defined as in Section \ref{plum-blossom-multiplication}, so that the all integers becomes a commutative semigroup under this operation. The plum-blossom product method not only can be efficiently used to mental calculation but also can be applied to design new multiplier, see \cite{Zhu1}.

In the present paper, as the development of the plum-blossom product method, we propose the plum-blossom wedge product method, which may be more efficiently applied to multiplication and division of long integers, especially to mental arithmetic.

This article is organized as follows. Section \ref{introduction} is an introduction and Section \ref{222} is the basic formula of the rapid multiplication of two integers.
A brief introduction to the plum-blossom product method of multiplication and division of multi-digit numbers in Section \ref{plum-blossom-multiplication} and Section \ref{plum-blossom-division}, respectively.
In Section \ref{def-wedge-product} and Section \ref{prop-wedge-product},
the definition and some fundamental of plum-blossom wedge products are introduced respectively.
For every 1-digit multiplier, the plum-blossom wedge product table is presented in Section \ref{wedge-table}.
Next, for any long integers, we propose the plum-blossom wedge product method of multiplication in Section \ref{wedge-product-multiplication} and division in Section \ref{wedge-product-division}. At last, the conclusion is in Section \ref{conclusion}.

\section{The basic formula of the rapid multiplication of two integers}\label{222}

As preparation, we present the basic formula of the rapid multiplication of two integers, which was in effect used by some authors, see for example, \cite{Kavita,Shi}. But here is the most specific and most general form of the formula.
\par
Let's begin with a notion and its notation. For any two groups of numbers $a_1,a_2,\ldots,a_n$ and $b_1,b_2,\ldots, b_n$, the \emph{cross product sum} is defined as
$$\begin{bmatrix} a_1 & a_2 & \cdots & a_n \\ b_1 & b_2 & \cdots & b_n  \end{bmatrix}=\sum_{i=1}^n a_i\times b_{n+1-i}=a_1\times b_n+a_2\times b_{n-1}+\cdots+a_n\times b_1.$$
For example, the cross product sum of $1,2,3$ and $4,5,6$ is $$\begin{bmatrix} 1 & 2 &3 \\4 & 5  & 6 \end{bmatrix}=1\times 6 +2\times 5 +3\times 4=28.$$
Evidently, the cross product sum of two numbers is exactly their product, that is,
$$\begin{bmatrix} a\\b\end{bmatrix}=a\times b.$$
\par
For example, to calculate the product $123\times 456$, we can use the cross product sums.
If we momentarily ignore the carry of all involved cross product sums, then the units of this  product is equal to $3\times 6$, that is $\begin{bmatrix} 3\\6 \end{bmatrix}$;
the tens is $2\times 6+3\times 5$, that is $\begin{bmatrix} 2 & 3\\5 & 6 \end{bmatrix}$;
the hundreds is $1\times 6+2\times 5+3\times 4$, that is $\begin{bmatrix} 1 & 2 & 3\\4 & 5 & 6 \end{bmatrix}$;
the thousands is $1\times 5+2\times 4$, that is $\begin{bmatrix} 1 & 2 \\4 & 5 \end{bmatrix}$;
the ten-thousands is $1\times 4$, that is $\begin{bmatrix} 1 \\4\end{bmatrix}$.
Therefore, $$123\times 456=(\begin{bmatrix} 1 \\4\end{bmatrix},\begin{bmatrix} 1 & 2 \\4 & 5 \end{bmatrix},\begin{bmatrix} 1 & 2 & 3\\4 & 5 & 6 \end{bmatrix},\begin{bmatrix} 2 & 3\\5 & 6 \end{bmatrix},\begin{bmatrix} 3\\6 \end{bmatrix}).$$
The final result of the product can be obtained by calculating each cross product sum and carrying out the rounding process:
$$123\times 456=(4,13,28,27,18)=(5,5,t,8,8)=56\,088,$$
where $t$ represents the number $10$, which causes a further carry.
\par

Generally, we have
\begin{lem}[The Basic Formula of the Rapid Multiplication ]\label{lem:chengjijibengongshi-1}\label{basic-formula}
For $m\geq n$, the product of an $m$-digit number $(a_1,a_2,\ldots, a_m)$ and an $n$-digit number $(b_1,b_2,\ldots,b_n)$ equals to
\begin{gather*}
(a_1,a_2,\ldots, a_m)\times (b_1,b_2,\ldots,b_n)\\
=(\begin{bmatrix} a_1 \\b_1\end{bmatrix},\begin{bmatrix} a_1 & a_2 \\b_1 & b_2 \end{bmatrix},\cdots, \begin{bmatrix} a_1 & a_2 & \ldots & a_n\\b_1 & b_2 & \ldots & b_n\end{bmatrix},
\begin{bmatrix} a_2 & a_3 & \ldots & a_{n+1}\\b_1 & b_2 & \ldots & b_n\end{bmatrix},
\ldots,
\\
\begin{bmatrix} a_{m-n+1} & a_{m-n+2} & \ldots & a_m\\b_1 & b_2 & \ldots & b_n\end{bmatrix},
\begin{bmatrix} a_{m-n+2} & a_{m-n+3} & \ldots & a_m\\b_2 & b_3 & \ldots & b_n\end{bmatrix},
\cdots,
\begin{bmatrix} a_{m-1} & a_m \\b_{n-1} & b_n \end{bmatrix},\begin{bmatrix} a_m\\b_n \end{bmatrix}).
\end{gather*}
\end{lem}

For multi-digit numbers, we may use the segmented numbers. For example, for the six-digit number $123\,456$, the segmentations by length $2$ are $(12,34,56)$, and the segmentations by length $3$ are $(123,456)$. Notice that, the above basic formula (in Lemma~\ref{lem:chengjijibengongshi-1}) is also applicable for segmented numbers whose segment length is greater than $1$. Here is an example to use the basic formula of the rapid multiplication to segmented numbers:

\begin{exmp}\rm
Using Lemma~\ref{lem:chengjijibengongshi-1} to segmented numbers with segment length $2$, we can compute the following multiplication in mind:
\begin{gather*}
2\,976\times 2\,924=(29,76)\times (29,24)=
(\begin{bmatrix} 29  \\ 29 \end{bmatrix},\begin{bmatrix} 29 & 76 \\ 29 & 24 \end{bmatrix},\begin{bmatrix}  76 \\ 24 \end{bmatrix})\\
=(\begin{bmatrix} 29  \\ 29 \end{bmatrix},\begin{bmatrix} 29 \\ 76+24 \end{bmatrix},\begin{bmatrix}  76 \\ 24 \end{bmatrix})
=(\begin{bmatrix} 29  \\ 29 \end{bmatrix},\begin{bmatrix} 29 & 100 \end{bmatrix},\begin{bmatrix}  76 \\ 24 \end{bmatrix})
=(\begin{bmatrix} 29  \\ 29+1 \end{bmatrix},0,\begin{bmatrix}  76 \\ 24 \end{bmatrix})\\
=(\begin{bmatrix} 29  \\ 30 \end{bmatrix},0,\begin{bmatrix}  76 \\ 24 \end{bmatrix})
=(\begin{bmatrix} 29  \\ 30 \end{bmatrix},0,\begin{bmatrix}  76 \\ 24 \end{bmatrix})
=(870,00,1824)=(8,70,18,24)=8\,701\,824.
\end{gather*}
\end{exmp}

\section{The plum-blossom product method of multiplication of multi-digit numbers}\label{plum-blossom-multiplication}

As in \cite{Zhu1,Zhu2,Zhu4}, we define plum-blossom product as below:

\begin{defn}
For two integers $a$ and $b$, the plum-blossom product of $a$ and $b$ is denoted as $a\clubsuit b$
and  is defined as the ones of the usual product of $a$ and $b$ if the ones is less than or equal to $3$ or otherwise the ones minus $10$.
\end{defn}

By the definition, the plum-blossom product of any two digits is less than or equal to $3$ and meanwhile it is greater than or equal to $-6$.
For example, since $3\times 7=21$ and $1\leqslant 3$, the plum-blossom product of $3$ and $7$ is $1$, that is, $3\clubsuit 7=1$. Since  $7\times 7=49$ and $9\ge 3$,  the plum-blossom product of $7$ and $7$ is $9-10=-1$, that is, $7\clubsuit 7=-1$. Similarly, $5\clubsuit 5=5-10=-5$ and $5\clubsuit 8=0$. From the usual multiplication table, one may easily deduce the following plum-product table:

\begin{table}[!htb]
\centering
 \begin{tabular} 
 {|c||c|c|c|c|c|c|c|c|c|c|}
                    \hline
                    $\clubsuit$  &  \rule{8pt}{0pt}$1$\rule{8pt}{0pt}&  \rule{8pt}{0pt}$2$\rule{8pt}{0pt}&  \rule{8pt}{0pt}$3$\rule{8pt}{0pt}&  \rule{8pt}{0pt}$4$\rule{8pt}{0pt}&  \rule{8pt}{0pt}$5$\rule{8pt}{0pt}&  \rule{8pt}{0pt}$6$\rule{8pt}{0pt}&  \rule{8pt}{0pt}$7$\rule{8pt}{0pt}&  \rule{8pt}{0pt}$8$\rule{8pt}{0pt}&  \rule{8pt}{0pt}$9$\rule{8pt}{0pt}\\
                                   \hline\hline
                   $1 $&$ 1 $&$ 2 $&$ 3 $&$ -6 $&$ -5 $&$ -4 $&$ -3 $&$ -2 $&$ -1 $\\
                   \hline
                   $2 $&$  $&$ -6 $&$ -4 $&$ -2 $&$ 0 $&$ 2 $&$ -6 $&$ -4 $&$ -2 $\\
                   \hline
                  $ 3 $&$  $&$  $&$ -1 $&$ 2 $&$ -5$ &$ -2 $&$ 1 $&$ -6 $&$ -3 $\\
                   \hline
                  $ 4 $&$  $&$  $&$  $&$ -4 $&$ 0 $&$ -6 $&$ -2 $&$ 2 $&$ -4 $\\
                   \hline
                   $5 $&$  $&$  $&$  $&$     $&$-5 $&$ 0 $&$ -5 $&$ 0 $&$ -5 $\\
                   \hline
                   $6 $&$  $&$  $&$  $&$  $&$  $&$ -4 $&$ 2 $&$ -2 $&$ -6 $\\
                   \hline
                   $7 $&$  $&$  $&$  $&$  $&$  $&$  $&$ -1  $&$ -4 $&$  3 $\\
                   \hline
                  $ 8 $&$  $&$  $&$  $&$  $&$  $&$  $&$  $&$ -6 $&$ 2 $\\
                   \hline
                   $9 $  &$  $&$  $&$  $&$  $&$  $&$  $&$  $&$  $&$ 1 $\\
                    \hline
 \end{tabular}
 \caption{The plum-blossom product table}
\label{tab:jiandaoji}
\end{table}

Since the plum-blossom product satisfies the commutativity, i.e. $a\clubsuit b=b\clubsuit a$, if all blanks in above table are filled in completely, it can be seen that the whole table is completely symmetric with respect to its main diagonal. It is due to this symmetry that we have left half of our table blank. The advantage of doing so, of course, is that it makes the table more cleaner. It is seen that the third of $45$ plum-blossom products is negative, which ensures that there is a good chance that the pluses and minuses will balance out when calculating the sum of the cross plum-blossom products and the carry from the next position.

\begin{thm}\label{thm-clubsuit}
Given two digits $a\leq b$, suppose that $a\times b=(J(a\clubsuit b),a\clubsuit b)$. Then
\begin{equation*}
J(a\clubsuit b)=
\begin{cases}
a & \textrm{if}\; a=1 \;\textrm{or} \; b=9, \textrm{and} \; b-a\ge 3;\\
a & \textrm{if}\; b-a\ge 5;\\
a-2 &\textrm{if}\; 3\le a\le b\le 7 \;\textrm{and}\; b-a\le 1;\\
a-1 & \;\textrm{otherwise}.
\end{cases}
\end{equation*}
\end{thm}

Let $J(a\clubsuit b)=\min(a,b)+\delta(a,b)$. Then the above theorem tells us that $\delta(a,b)=0,-1$ or $-2$.

Theorem \ref{thm-clubsuit} incorporating into the basic formula of the rapid multiplication provides a novel multiplication method, which we called the plum-blossom product method. A great amount of examples shows that in most cases the absolute value of the sum of the cross plum-blossom products and the carry from the next position would not be too large, which provides a great convenience for mental multiplication. Let's  see two examples.

\begin{exmp}\rm
Compute $386\times 47. $
\end{exmp}
\par
\solution \hskip 8pt
First by the basic formula of the rapid multiplication, we have $$\begin{bmatrix} 386\\47\end{bmatrix}=(\begin{bmatrix}3\\4 \end{bmatrix},\begin{bmatrix}3 & 8 \\ 4 &7 \end{bmatrix},\begin{bmatrix} 8 & 6 \\ 4 &7 \end{bmatrix}, \begin{bmatrix}6\\7 \end{bmatrix})$$
Then by the plum-blossom product method, we obtain
\begin{align*}
  & 386\times 47\\
=&(3\times 4+3+4-2,3\clubsuit 7+8\clubsuit 4+4+7-2,8\clubsuit 7+6\clubsuit 4+6-1-1,6\clubsuit 7)\\
=&(17,12,-6,2)\\
=&18142.
\end{align*}
\endsolution

\begin{exmp}\rm
Compute $456\times 789. $
\end{exmp}
\par
\solution \hskip 8pt
First by the basic formula of the rapid multiplication, we have
$$\begin{bmatrix} 456\\789\end{bmatrix}
=(\begin{bmatrix}4\\7 \end{bmatrix},\begin{bmatrix}4 & 5 \\ 7 & 8 \end{bmatrix},\begin{bmatrix}4 & 5 & 6\\ 7 & 8 & 9\end{bmatrix},
\begin{bmatrix} 5 & 6 \\ 8 & 9 \end{bmatrix}, \begin{bmatrix}8 \\ 9 \end{bmatrix})$$
Then by the plum-blossom product method, we obtain
\begin{align*}
  & 456\times 789\\
=&(4\times 7+4+5-2,4\clubsuit 8+5\clubsuit 7+4+5+6-3,4\clubsuit 9+5\clubsuit 8+6\clubsuit 7+5+6-2+1,\\
 & 5\clubsuit 9+6\clubsuit 8,6\times 9)\\
=&(35,9,8,-2,4)\\
=&359\,784.
\end{align*}
\endsolution

The plum-blossom products have some fundamental properties, which we present as below.

\begin{lem}\label{lem-clubsuit-property}
For any integers $a$, $b$, $c$, the following statements hold:
\begin{enumerate}[(1)]
\item $a\clubsuit b=b\clubsuit a$;
\item $(a\clubsuit b)\clubsuit c=a\clubsuit (b\clubsuit c)$;
\item $(a\times b)\clubsuit c=a\clubsuit (b\times c)$;
\item $a\clubsuit b=a\clubsuit (b-10)=(a-10)\clubsuit b=(a-10)\clubsuit (b-10)$;
\item if $a$ is an even number, then $a\clubsuit b=a\clubsuit (b-5)$;
\item if the one of $a$ and $b$ is an odd number and the other is even number, then $a\clubsuit b=(a-5)\clubsuit (b-5)$.
\end{enumerate}
\end{lem}

\section{The plum-blossom product method for division of multi-digit integers}\label{plum-blossom-division}

Since division is the inverse operation of multiplication, from the plum-blossom product method for multiplication,
we deduce the plum-blossom product method for division of  multi-digit inters. But this process of yielding is not direct, since it requires certain techniques which we present below.

Let us consider the division $a\div b$ of multi-digits $a=(a_1,a_2,\ldots, a_s)$ and $b=(b_1, b_2, \ldots, b_t)$.
Suppose that the integral part of $a\div b$ is $c=(c_1, c_2, \ldots, c_{s-t+1})$ where it is allowable that $c_1=0$.
For convenience, set $c_n=0$ for all $n>s-t+1$.
Similar with the usual vertical division method, starting with the highest digit of the dividend, we obtain the quotient $c$ step by step. It is different from the usual method, we get the partial remainders by subtracting the partial plum-blossom product of $b\times c$ at each digit position step by step.
At each step, we do the subtraction in two steps. The one step is subtracting the partial plum-blossom product not involving in the newest quotient $c_n$, that is,
\[PP_n^{(0)}=b_2\clubsuit c_{n-1}+b_3\clubsuit c_{n-2}+\cdots+J(b_3\clubsuit c_{n-1})+J(b_4\clubsuit c_{n-2})+\cdots,\]
and the second step, after determine the n'th digit $c_n$ of the quotient, is subtracting the partial plum-blossom product involving in the newest quotient $c_n$, that is,
$$PP_n^{(1)}=b_1\times c_n+J(b_1\clubsuit c_{n-1}),$$
where $$b_1\times c_n=b_1\clubsuit c_n+J(b_1\clubsuit c_n)\times 10.$$
Thus we have
\[PP_n=PP_n^{(0)}+PP_n^{(1)},\quad \text{for}\quad n=1,2,\ldots, s.\]
It is evident that $PP_1^{(0)}=0$. It is also easily seen that $PP_n^{(1)}=0$ for all $n>s-t+1.$
The partial remainder, denoted by $r_m$, is defined as follows:
\[r_1=a_1-PP_1\]
\[r_n=(r_{n-1},a_n)-PP_n \quad \text{for each}\quad n>1.\]
It is obvious that the last partial remainder $r_s$ is exactly the final remainder of the division $a\div b$, that is,
\[a\div b=c\cdots\cdots r_s.\]
The above  algorithm of division is called the \emph{plum-blossom product method of division}, or the \emph{plum-blossom product division}.

\begin{exmp}\rm
Compute the division $56789\div 369$ with the plum-blossom product method.
\end{exmp}
\par
{\bf Solution\quad}
Since $s-t+1=5-3+1=3$, we can suppose that the quotient is $c=(c_1,c_2,c_3)$.
\begin{enumerate}[$\bullet$]
\item {\it Step 1.} $a_1=5, b_1=3, c_1=1$ so that $PP_1=PP_1^{(1)}=1\times 3+J(1\clubsuit6)=3+1=4$ and that $r_1=5-4=1$.

\item {\it Step 2.} $a_2=6$ so that $(r_1,a_2)=16$; $PP_2^{(0)}=6\clubsuit 1+J(9\clubsuit 1)=-4+1=-3$ and $(r_1,a_2)-PP_2^{(0)}=16-(-3)=19$; $c_2=5$ so that
$PP_2^{(1)}=3\times 5+J(6\clubsuit 5)=15+3=18$ and that $r_2=19-18=1$.

\item {\it Step 3.} $a_3=7$ so that $(r_2,a_3)=17$; $PP_3^{(0)}=6\clubsuit 5+9\clubsuit 1+J(9\clubsuit 5)=0-1+4=4$ and $(r_2,a_3)-PP_3^{(0)}=17-4=13$; $c_3=3$ so that
$PP_3^{(1)}=3\times 3+J(6\clubsuit 3)=9+2=11$ and that $r_3=13-11=2$.

\item {\it Step 4.} $a_4=8$ so that $(r_3,a_4)=28$; $PP_4=PP_4^{(0)}=6\clubsuit 3+9\clubsuit 5+J(9\clubsuit 3)=-2-5+3=-4$ and $r_4=(r_3,a_4)-PP_4=28-(-4)=32$.

\item {\it Step 5.} $a_5=9$ so that $(r_4,a_5)=329$; $PP_5=PP_5^{(0)}=9\clubsuit 3=-3$ and $r_5=(r_4,a_5)-PP_5=329-(-3)=332$.

\end{enumerate}
Therefore, we obtain
\[56789\div 369=153\cdots\cdots332.\]
The above operation process can be presented as in the following vertical division:
$$
\begin{array}{lr}
&153\enspace\enspace \\
369 \!\!\!\!\! & \overline{)\enspace 56789} \\
& \underline{4} \enspace\enspace\enspace\enspace \\
& 16 \enspace\enspace\enspace \\
& \underline{\enspace \overline{3}}  \enspace\enspace\enspace \\
& 19 \enspace\enspace\enspace \\
& \underline{18} \enspace\enspace\enspace\\
& 17 \enspace\enspace \\
& \underline{\enspace 4}  \enspace\enspace \\
& 13 \enspace\enspace \\
& \underline{11} \enspace\enspace \\
& 28 \enspace \\
& \underline{\enspace\overline{4}}  \enspace\\
& 329 \\
&\underline{\enspace\enspace\overline{3}}\\
& 332 \end{array}
$$

\begin{exmp}\rm
Compute the division $56789\div 369$ with the plum-blossom product method so that the result is account to two decimal places.
\end{exmp}
\par
{\bf Solution\quad}
Since required being correct to two decimal places, we multiply the dividend by 100 to obtain 5678900 so that the question is turned to computing the division $5678900\div 369$.
We use the vertical division below so that we obtain that
\[5678900\div 369=15389\cdots\cdots359,\]
divided by 100, which yields that
\[56789\div 369=153.89\cdots\cdots3.59.\]

$$
\begin{array}{lr}
&15389\enspace\enspace \\
369 \!\!\!\!\! & \overline{)\enspace 5678900} \\
& \underline{4} \enspace\enspace\enspace\enspace\enspace\enspace \\
& 16 \enspace\enspace\enspace\enspace\enspace \\
& \underline{\enspace \overline{3}}  \enspace\enspace\enspace\enspace\enspace \\
& 19 \enspace\enspace\enspace\enspace\enspace \\
& \underline{18} \enspace\enspace\enspace\enspace\enspace\\
& 17 \enspace\enspace\enspace\enspace \\
& \underline{\enspace 4}  \enspace\enspace\enspace\enspace \\
& 13 \enspace\enspace\enspace\enspace \\
& \underline{11} \enspace\enspace\enspace\enspace \\
& 28 \enspace\enspace\enspace \\
& \underline{\enspace\overline{4}}  \enspace\enspace\enspace\\
& 32\enspace\enspace\enspace \\
& \underline{29}\enspace\enspace\enspace \\
& 39\enspace\enspace \\
& \underline{\enspace2}\enspace\enspace \\
& 37\enspace\enspace \\
& \underline{33}\enspace\enspace \\
& 40\enspace \\
& \underline{\enspace4}\enspace \\
& 360\\
& \underline{\enspace\enspace1}\\
& 359 \end{array}
$$

\begin{exmp}\rm
Compute the division $2728018\div 3456$ with the plum-blossom product method.
\end{exmp}
\par
{\bf Solution\quad}
Noting that $s-t+1=7-4+1=4$,  we suppose that the quotient is $c=(c_1,c_2,c_3,c_4)$. Since $a_1=2$ and $b_1=3$, $a_1<b_1$ so that $c_1=0$. The whole calculation process is as in the following vertical division, which gives the answer:
\[2728018\div 3456=789\cdots\cdots 1234.\]
$$
\begin{array}{lr}
&789\enspace\enspace\enspace \\
3456 \!\!\!\!\! & \overline{)\enspace 2728018} \\
& \underline{24} \enspace\enspace\enspace\enspace\enspace \\
& 32 \enspace\enspace\enspace\enspace \\
& \underline{\enspace 2}  \enspace\enspace\enspace\enspace \\
& 30 \enspace\enspace\enspace\enspace \\
& \underline{27} \enspace\enspace\enspace\enspace\\
& 38 \enspace\enspace\enspace \\
& \underline{\enspace 5}  \enspace\enspace\enspace \\
& 33 \enspace\enspace\enspace \\
& \underline{31} \enspace\enspace\enspace \\
& 20 \enspace\enspace \\
& \underline{\enspace 8}  \enspace\enspace\\
& 121\enspace \\
& \underline{\enspace\enspace\overline{1}}\enspace \\
& 1228 \\
&\underline{\enspace\enspace\enspace\overline{6}}\\
& 1234 \end{array}
$$

\section{The definition of plum-blossom wedge products}\label{def-wedge-product}

Now we introduce the main concept of this paper, which provides a new method for multiplication and division of multi-digit numbers.

\begin{defn}
Given digits $a$, $b$ and $c$, the plum-blossom wedge product of $(a,b)$ and $c$ is denoted as $(a,b)\bowtie c$ and defined by
$$(a,b)\bowtie c=a\clubsuit c+J(b\clubsuit c).$$
\end{defn}

For example, we have
$$35\bowtie 7=3\clubsuit 7+J(5\clubsuit 7)=1+4=5,$$
$$46\bowtie 8=4\clubsuit 8+J(6\clubsuit 8)=2+5=7,$$
$$59\bowtie 7=5\clubsuit 7+J(9\clubsuit 7)=-5+6=1,$$
$$46\bowtie 9=4\clubsuit 9+J(6\clubsuit 9)=-4+6=2,$$
$$74\bowtie 2=7\clubsuit 2+J(4\clubsuit 2)=-6+1=-5.$$

\section{The properties of plum-blossom wedge products}\label{prop-wedge-product}

The plum-blossom wedge products have some fine features which we discuss in this section.

\begin{prop}
For any digits $a$, $b$ and $c$, we have
\begin{enumerate}[(1)]
\item $-6\le (a,b)\bowtie c\le 11$;
\item if $c\ne 9$, then $(a,b)\bowtie c\le 9$.
\end{enumerate}
\end{prop}
In fact, $78\bowtie 9=10$ and $79\bowtie 9=11$. Accept the both cases, all plum-blossom wedge products are less than 9.

By Lemma \ref{lem-clubsuit-property}, one may easily deduce the following proposition for plum-blossom wedge products.

\begin{prop}
For any digits $a$, $b$ and $c$, we have
\begin{enumerate}[(1)]
\item if $c$ is an even number and $a\le 4$, then $(a+5,b)\bowtie c=(a,b)\bowtie c$;
\item if $c$ is an even number and $b\le 4$, then $(a,b+5)\bowtie c=(a,b)\bowtie c+\frac{c}{2}$;
\item for fixed $a$ and $c$,  $(a+5,b)\bowtie c$ is monotone increasing, that is, $(a,b)\bowtie c\le (a,b+1)\bowtie c$;
\item for fixed $a$ and $c$,  $(a+1,b)\bowtie c=(a,b)\bowtie c+c$ or $(a+1,b)\bowtie c=(a,b)\bowtie c+c-10$.
\end{enumerate}
\end{prop}

The main results of this section is the following two theorems.

\begin{thm}
For any digits $a$ and $b$, suppose that $a\times b=(c,d)$ where $0\le c,d\le 9$. Then
\begin{enumerate}[(1)]
\item if $d\le 3$, then $(a, a)\bowtie b=c+d$;
\item if $d>3$, then $(a, a)\bowtie b=c+d-9$.
\end{enumerate}
\end{thm}

\begin{proof}
Suppose that $d\le 3$. Then $a\clubsuit b=d$ and $J(a\clubsuit b)=c$. It follows that $(a, a)\bowtie b=c+d$, that is, (1) holds.\par
suppose that $d>3$. Then $a\clubsuit b=d-10$ and $J(a\clubsuit b)=c+1$. It follows that $(a, a)\bowtie b=(c+1)+(d-10)=c+d-9$, which means that (2) follows.
\end{proof}

For example, since $3\times 7=21$ and $1\le 3$, $33\bowtie 7=2+1=3$. Since $3\times 8=24$ and $4>3$, $33\bowtie 8=2+4-9=-3$.

Denote $a^+=a+1$ and recall that $J(a\clubsuit b)=\min(a,b)+\delta(a,b)$. Then we have the following important result for plum-bossom wedge products.
\begin{thm}
Let $a$, $b$ and $c$ be digits with that $b\ge c$. Then the following statements hold.
\begin{enumerate}[(1)]
\item if $a\ge c$ and $\delta(a^+,c)=\delta(a,c)$, or if $a<c$ and $\delta(a^+,c)\ne \delta(a,c)$. Then $$(a, b)\bowtie c=a^+\clubsuit c+\delta(b,c);$$
\item if $a\ge c$ and $\delta(a^+,c)\ne \delta(a,c)$, or if $a<c$ and $\delta(a^+,c)=\delta(a,c)$. Then $$(a, b)\bowtie c=(a^+\clubsuit c+10)+\delta(b,c).$$
\end{enumerate}
\end{thm}

\begin{proof}
Let $a$, $b$ and $c$ be digits with that $b\ge c$. Then $J(b\clubsuit c)=c+\delta(b,c)$. Note that $a\clubsuit c=a\times c-J(a\clubsuit c)\times 10$.
It follows that
$$(a,b)\bowtie c=a\clubsuit c+J(b\clubsuit c)=a\times c-J(a\clubsuit c)\times 10+c+\delta(b,c)=(a+1)\times c-J(a\clubsuit c)\times 10+\delta(b,c).$$
\par
Suppose that $a\ge c$. Then $a^+>c$. Thus $J(a\clubsuit c)=c+\delta(a,c)$ and $J(a^+\clubsuit c)=c+\delta(a^+,c)$.
If $\delta(a^+,c)=\delta(a,c)$, then
$a^+\clubsuit c=a^+\times c-J(a^+\clubsuit c)\times 10=a^+\times c-J(a\clubsuit c)\times 10$.
It follows that
$$(a,b)\bowtie c=a^+\clubsuit c+\delta(b,c).$$
If $\delta(a^+,c)\ne \delta(a,c)$, then $\delta(a^+,c)=\delta(a,c)+1$ and thus
$$a^+\clubsuit c=(a+1)\times c-J(a^+\clubsuit c)\times 10=(a+1)\times c-(J(a\clubsuit c)+1)\times 10=(a+1)\times c-J(a\clubsuit c)\times 10-10. $$
This means that
$$a^+\clubsuit c+10=(a+1)\times c-J(a\clubsuit c)\times 10. $$
It follows that
$$(a,b)\bowtie c=(a^+\clubsuit c+10)+\delta(b,c).$$
\par
Suppose that $a<c$. The similar argument shows that if $\delta(a^+,c)\ne \delta(a,c)$ then $$(a, b)\bowtie c=a^+\clubsuit c+\delta(b,c),$$
and if $a<c$ and $\delta(a^+,c)=\delta(a,c)$ then
$$(a, b)\bowtie c=(a^+\clubsuit c+10)+\delta(b,c),$$
which completes the proof.
\end{proof}

\section{The plum-blossom wedge product tables}\label{wedge-table}

For fixed digit $c$ and for all two-digit positive integers $(a,b)$, one can make the tables of  plum-blossom wedge products $(a,b)\bowtie c$ as follows, where the first column represents the value of $a$ and the first row represents the value of $b$.

\begin{table}[!htb]
\centering
 \begin{tabular} 
 {|c||c|c|c|c|c|c|c|c|c|c|c|}
                    \hline
                    $\bowtie(c=2)$ & \rule{8pt}{0pt}$0$\rule{8pt}{0pt} &  \rule{8pt}{0pt}$1$\rule{8pt}{0pt}&  \rule{8pt}{0pt}$2$\rule{8pt}{0pt}&  \rule{8pt}{0pt}$3$\rule{8pt}{0pt}&  \rule{8pt}{0pt}$4$\rule{8pt}{0pt}&  \rule{8pt}{0pt}$5$\rule{8pt}{0pt}&  \rule{8pt}{0pt}$6$\rule{8pt}{0pt}&  \rule{8pt}{0pt}$7$\rule{8pt}{0pt}&  \rule{8pt}{0pt}$8$\rule{8pt}{0pt}&  \rule{8pt}{0pt}$9$\rule{8pt}{0pt}\\
                                   \hline\hline
                   $0 $&$ 0 $&$ 0 $&$ 1 $&$1$&$1$&$1$&$1$&$2$&$ 2$& $2$\\
                   \hline
                   $1 $&$ 2 $&$ 2 $&$ 3 $&$ 3$&$3$&$3$&$ 3$&$ 4 $&$ 4 $& $4$\\
                   \hline
                   $2 $&$ -6 $&$-6$&$-5$&$-5$&$-5$&$-5$&$-5$&$-4$&$-4$& $-4$\\
                   \hline
                  $ 3$&$ -4 $&$-4$&$-3$&$-3$&$-3$&$-3$&$-3$&$-2$&$-2$& $-2$\\
                   \hline
                  $ 4 $&$ -2 $&$-2$&$-1$&$-1$&$-1$&$-1$&$-1$&$0$&$0$& $0$\\
                   \hline
                   $5$&$ 0 $&$ 0 $&$ 1 $&$1$&$1$&$1$&$1$&$2$&$ 2$& $2$\\
                   \hline
                   $6$&$ 2 $&$ 2 $&$ 3 $&$ 3$&$3$&$3$&$ 3$&$ 4 $&$ 4 $& $4$\\
                   \hline
                   $7 $&$ -6 $&$-6$&$-5$&$-5$&$-5$&$-5$&$-5$&$-4$&$-4$& $-4$\\
                   \hline
                  $ 8$&$ -4 $&$-4$&$-3$&$-3$&$-3$&$-3$&$-3$&$-2$&$-2$& $-2$\\
                   \hline
                   $9 $&$ -2 $&$-2$&$-1$&$-1$&$-1$&$-1$&$-1$&$0$&$0$& $0$\\
                    \hline
 \end{tabular}
 \caption{The plum-blossom wedge product table for $c=2$}
\label{tab:jiandaoji}
\end{table}

From the table for $c=2$ we see that
for any $a$, $(a,0)\bowtie 2=(a,1)\bowtie 2$, $(a,2)\bowtie 2=(a,3)\bowtie 2=(a,4)\bowtie 2=(a,5)\bowtie 2=(a,6)\bowtie 2$,
$(a,7)\bowtie 2=(a,8)\bowtie 2=(a,9)\bowtie 2$.

\begin{table}[!htb]
\centering
 \begin{tabular} 
 {|c||c|c|c|c|c|c|c|c|c|c|c|}
                    \hline
                    $\bowtie(c=3)$ & \rule{8pt}{0pt}$0$\rule{8pt}{0pt} &  \rule{8pt}{0pt}$1$\rule{8pt}{0pt}&  \rule{8pt}{0pt}$2$\rule{8pt}{0pt}&  \rule{8pt}{0pt}$3$\rule{8pt}{0pt}&  \rule{8pt}{0pt}$4$\rule{8pt}{0pt}&  \rule{8pt}{0pt}$5$\rule{8pt}{0pt}&  \rule{8pt}{0pt}$6$\rule{8pt}{0pt}&  \rule{8pt}{0pt}$7$\rule{8pt}{0pt}&  \rule{8pt}{0pt}$8$\rule{8pt}{0pt}&  \rule{8pt}{0pt}$9$\rule{8pt}{0pt}\\
                                   \hline\hline
                   $0$ & $0$ & $0$ & $1$ & $1$ & $1$ & $2$ & $2$ & $2$ & $3$ & $3$ \\
                   \hline
                   $1$ & $3$ & $3$ & $4$ & $4$ & $4$ & $5$ & $5$ & $5$ & $6$ & $6$ \\
                   \hline
                   $2$ & $-4$ & $-4$ & $-3$ & $-3$ & $-3$ & $-2$ & $-2$ & $-2$ & $-1$ & $-1$ \\
                   \hline
                   $3$ & $-1$ & $-1$ & $0$ & $0$ & $0$ & $1$ & $1$ & $1$& $2$ & $2$ \\
                   \hline
                   $4$ & $2$ & $2$ & $3$ & $3$ & $3$ & $4$ & $4$ & $4$ & $5$ & $5$ \\
                   \hline
                   $5$ & $-5$ & $-5$ & $-4$ & $-4$ & $-4$ & $-3$ & $-3$ & $-3$ & $-2$ & $-2$ \\
                   \hline
                   $6$ & $-2$ & $-2$ & $-1$ & $-1$ & $-1$ & $0$ & $0$ & $0$ & $1$ & $1$ \\
                   \hline
                   $7$ & $1$ & $1$ & $2$ & $2$ & $2$ & $3$ & $3$ & $3$ & $4$ & $4$ \\
                   \hline
                   $8$ & $-6$ & $-6$ & $-5$ & $-5$ & $-5$ & $-4$ & $-4$ & $-4$ & $-3$ & $-3$ \\
                   \hline
                   $9$ & $-3$ & $-3$ & $-2$ & $-2$ & $-2$ & $-1$ & $-1$ & $-1$ & $0$ & $0$ \\
                    \hline
 \end{tabular}
 \caption{The plum-blossom wedge product table for $c=3$}
\label{tab:jiandaoji}
\end{table}

From the table for $c=3$ we see the following statements hold.
\begin{enumerate}[(1)]
\item for any $a,b$ with that $b\le 6$,  $(a,b+3)\bowtie 3=(a,b)\bowtie 3+1$;
\item for any $a,b$ with that $a\le 6$, $(a+3,b)\bowtie 3=(a,b)\bowtie 3-1$;
\item for any $a,k$, $(a,3k-1)\bowtie 3=(a,3k)\bowtie 3=(a,3k+1)\bowtie 3$ where $3k\pm 1$ is digits.
\end{enumerate}

\begin{table}[!htb]
\centering
 \begin{tabular} 
 {|c||c|c|c|c|c|c|c|c|c|c|c|}
                    \hline
                    $\bowtie(c=4)$ & \rule{8pt}{0pt}$0$\rule{8pt}{0pt} &  \rule{8pt}{0pt}$1$\rule{8pt}{0pt}&  \rule{8pt}{0pt}$2$\rule{8pt}{0pt}&  \rule{8pt}{0pt}$3$\rule{8pt}{0pt}&  \rule{8pt}{0pt}$4$\rule{8pt}{0pt}&  \rule{8pt}{0pt}$5$\rule{8pt}{0pt}&  \rule{8pt}{0pt}$6$\rule{8pt}{0pt}&  \rule{8pt}{0pt}$7$\rule{8pt}{0pt}&  \rule{8pt}{0pt}$8$\rule{8pt}{0pt}&  \rule{8pt}{0pt}$9$\rule{8pt}{0pt}\\
                                   \hline\hline
                   $0$ & $0$ & $1$ & $1$ & $1$ & $2$ & $2$ & $3$ & $3$ & $3$ & $4$ \\
                   \hline
                   $1$ & $-6$ & $-5$ & $-5$ & $-5$ & $-4$ & $-4$ & $-3$ & $-3$ & $-3$ & $-2$ \\
                   \hline
                   $2$ & $-2$ & $-1$ & $-1$ & $-1$& $0$ & $0$ & $1$ & $1$ & $1$ & $2$ \\
                   \hline
                   $3$ & $2$ & $3$ & $3$ & $3$ & $4$ & $4$ & $5$ & $5$ & $5$ & $6$ \\
                   \hline
                   $4$ & $-4$ & $-3$ & $-3$ & $-3$ & $-2$ & $-2$ & $-1$ & $-1$ & $-1$ & $0$ \\
                   \hline
                   $5$ & $0$ & $1$ & $1$ & $1$ & $2$ & $2$ & $3$ & $3$ & $3$ & $4$ \\
                   \hline
                   $6$ & $-6$ & $-5$ & $-5$ & $-5$ & $-4$ & $-4$ & $-3$ & $-3$ & $-3$ & $-2$ \\
                   \hline
                   $7$ & $-2$ & $-1$ & $-1$ & $-1$& $0$ & $0$ & $1$ & $1$ & $1$ & $2$ \\
                   \hline
                   $8$ & $2$ & $3$ & $3$ & $3$ & $4$ & $4$ & $5$ & $5$ & $5$ & $6$ \\
                   \hline
                   $9$ & $-4$ & $-3$ & $-3$ & $-3$ & $-2$ & $-2$ & $-1$ & $-1$ & $-1$ & $0$ \\
                    \hline
 \end{tabular}
 \caption{The plum-blossom wedge product table for $c=4$}
\label{tab:jiandaoji}
\end{table}

From the table for $c=4$ we see that $(a,b)\bowtie 4=(a,b')\bowtie 4$ whenever $b$ and $b'$ are in the same sets $A$, $B$ or $C$
where $A=\{1,2,3\}$, $B=\{4,5\}$ and $C=\{6,7,8\}$.

\begin{table}[!htb]
\centering
 \begin{tabular} 
 {|c||c|c|c|c|c|c|c|c|c|c|c|}
                    \hline
                    $\bowtie(c=5)$ & \rule{8pt}{0pt}$0$\rule{8pt}{0pt} &  \rule{8pt}{0pt}$1$\rule{8pt}{0pt}&  \rule{8pt}{0pt}$2$\rule{8pt}{0pt}&  \rule{8pt}{0pt}$3$\rule{8pt}{0pt}&  \rule{8pt}{0pt}$4$\rule{8pt}{0pt}&  \rule{8pt}{0pt}$5$\rule{8pt}{0pt}&  \rule{8pt}{0pt}$6$\rule{8pt}{0pt}&  \rule{8pt}{0pt}$7$\rule{8pt}{0pt}&  \rule{8pt}{0pt}$8$\rule{8pt}{0pt}&  \rule{8pt}{0pt}$9$\rule{8pt}{0pt}\\
                                   \hline\hline
                   $0$ & $0$ & $1$ & $1$ & $2$ & $2$ & $3$ & $3$ & $4$ & $4$ & $5$ \\
                   \hline
                   $1$ & $-5$ & $-4$ & $-4$ & $-3$ & $-3$ & $-2$ & $-2$ & $-1$ & $-1$ & $0$ \\
                   \hline
                   $2$  & $0$ & $1$ & $1$ & $2$ & $2$ & $3$ & $3$ & $4$ & $4$ & $5$ \\
                   \hline
                   $3$ & $-5$ & $-4$ & $-4$ & $-3$ & $-3$ & $-2$ & $-2$ & $-1$ & $-1$ & $0$ \\
                   \hline
                   $4$  & $0$ & $1$ & $1$ & $2$ & $2$ & $3$ & $3$ & $4$ & $4$ & $5$ \\
                   \hline
                   $5$ & $-5$ & $-4$ & $-4$ & $-3$ & $-3$ & $-2$ & $-2$ & $-1$ & $-1$ & $0$ \\
                   \hline
                   $6$  & $0$ & $1$ & $1$ & $2$ & $2$ & $3$ & $3$ & $4$ & $4$ & $5$ \\
                   \hline
                   $7$ & $-5$ & $-4$ & $-4$ & $-3$ & $-3$ & $-2$ & $-2$ & $-1$ & $-1$ & $0$ \\
                   \hline
                   $8$  & $0$ & $1$ & $1$ & $2$ & $2$ & $3$ & $3$ & $4$ & $4$ & $5$ \\
                   \hline
                   $9$ & $-5$ & $-4$ & $-4$ & $-3$ & $-3$ & $-2$ & $-2$ & $-1$ & $-1$ & $0$ \\
                    \hline
 \end{tabular}
 \caption{The plum-blossom wedge product table for $c=5$}
\label{tab:jiandaoji}
\end{table}

From the table for $c=5$ we see that for all $a\le 7$,
$$(a,b)\bowtie 5=(a+2,b)\bowtie 5.$$

\begin{table}[!htb]
\centering
 \begin{tabular} 
 {|c||c|c|c|c|c|c|c|c|c|c|c|}
                    \hline
                    $\bowtie(c=6)$ & \rule{8pt}{0pt}$0$\rule{8pt}{0pt} &  \rule{8pt}{0pt}$1$\rule{8pt}{0pt}&  \rule{8pt}{0pt}$2$\rule{8pt}{0pt}&  \rule{8pt}{0pt}$3$\rule{8pt}{0pt}&  \rule{8pt}{0pt}$4$\rule{8pt}{0pt}&  \rule{8pt}{0pt}$5$\rule{8pt}{0pt}&  \rule{8pt}{0pt}$6$\rule{8pt}{0pt}&  \rule{8pt}{0pt}$7$\rule{8pt}{0pt}&  \rule{8pt}{0pt}$8$\rule{8pt}{0pt}&  \rule{8pt}{0pt}$9$\rule{8pt}{0pt}\\
                                   \hline\hline
                   $0$ & $0$ & $1$ & $1$ & $2$ & $3$ & $3$ & $4$ & $4$ & $5$ & $6$ \\
                   \hline
                   $1$ & $-4$ & $-3$ & $-3$ & $-2$ & $-1$ & $-1$ & $0$ & $0$ & $1$ & $2$ \\
                   \hline
                   $2$ & $2$ & $3$ & $3$ & $4$ & $5$ & $5$ & $6$ & $6$ & $7$ & $8$ \\
                   \hline
                   $3$ & $-2$ & $-1$ & $-1$ & $0$ & $1$ & $1$ & $2$ & $2$ & $3$ & $4$ \\
                   \hline
                   $4$ & $-6$ & $-5$ & $-5$ & $-4$ & $-3$ & $-3$ & $-2$ & $-2$ & $-1$ & $0$ \\
                   \hline
                   $5$ & $0$ & $1$ & $1$ & $2$ & $3$ & $3$ & $4$ & $4$ & $5$ & $6$ \\
                   \hline
                   $6$ & $-4$ & $-3$ & $-3$ & $-2$ & $-1$ & $-1$ & $0$ & $0$ & $1$ & $2$ \\
                   \hline
                   $7$ & $2$ & $3$ & $3$ & $4$ & $5$ & $5$ & $6$ & $6$ & $7$ & $8$ \\
                   \hline
                   $8$ & $-2$ & $-1$ & $-1$ & $0$ & $1$ & $1$ & $2$ & $2$ & $3$ & $4$ \\
                   \hline
                   $9$ & $-6$ & $-5$ & $-5$ & $-4$ & $-3$ & $-3$ & $-2$ & $-2$ & $-1$ & $0$ \\
                    \hline
 \end{tabular}
 \caption{The plum-blossom wedge product table for $c=6$}
\label{tab:jiandaoji}
\end{table}

From the table for $c=6$ we see the following statements hold.
\begin{enumerate}[(1)]
\item for any $a,b$ with that $b\le 6$,  $(a,b+3)\bowtie 6=(a,b)\bowtie 6+2$;
\item for any $a,b$ with that $a\le 6$, $(a+3,b)\bowtie 6=(a,b)\bowtie 6-2$;
\item for any $a\le 9$, $(a,1)\bowtie 6=(a,2)\bowtie 6$, $(a,4)\bowtie 6=(a,5)\bowtie 6$, $(a,6)\bowtie 6=(a,7)\bowtie 6$,
\end{enumerate}

\begin{table}[!htb]
\centering
 \begin{tabular} 
 {|c||c|c|c|c|c|c|c|c|c|c|c|}
                    \hline
                    $\bowtie(c=7)$ & \rule{8pt}{0pt}$0$\rule{8pt}{0pt} &  \rule{8pt}{0pt}$1$\rule{8pt}{0pt}&  \rule{8pt}{0pt}$2$\rule{8pt}{0pt}&  \rule{8pt}{0pt}$3$\rule{8pt}{0pt}&  \rule{8pt}{0pt}$4$\rule{8pt}{0pt}&  \rule{8pt}{0pt}$5$\rule{8pt}{0pt}&  \rule{8pt}{0pt}$6$\rule{8pt}{0pt}&  \rule{8pt}{0pt}$7$\rule{8pt}{0pt}&  \rule{8pt}{0pt}$8$\rule{8pt}{0pt}&  \rule{8pt}{0pt}$9$\rule{8pt}{0pt}\\
                                   \hline\hline
                   $0$ & $0$ & $1$ & $2$ & $2$ & $3$ & $4$ & $4$ & $5$ & $6$ & $6$ \\
                   \hline
                   $1$ & $-3$ & $-2$ & $-1$ & $-1$ & $0$ & $1$ & $1$ & $2$ & $3$ & $3$ \\
                   \hline
                   $2$ & $-6$ & $-5$ & $-4$ & $-4$ & $-3$ & $-2$ & $-2$ & $-1$ & $0$ & $0$ \\
                   \hline
                   $3$ & $1$ & $2$ & $3$ & $3$ & $4$ & $5$ & $5$ & $6$ & $7$ & $7$ \\
                   \hline
                   $4$ & $-2$ & $-1$ & $0$ & $0$ & $1$ & $2$ & $2$ & $3$ & $4$ & $4$ \\
                   \hline
                   $5$ & $-5$ & $-4$ & $-3$ & $-3$ & $-2$ & $-1$ & $-1$ & $0$ & $1$ & $1$ \\
                   \hline
                   $6$ & $2$ & $3$ & $4$ & $4$ & $5$ & $6$ & $6$ & $7$ & $8$ & $8$ \\
                   \hline
                   $7$ & $-1$ & $0$ & $1$ & $1$ & $2$ & $3$ & $3$ & $4$ & $5$ & $5$ \\
                   \hline
                   $8$ & $-4$ & $-3$ & $-2$ & $-2$ & $-1$ & $0$ & $0$ & $1$ & $2$ & $2$ \\
                   \hline
                   $9$ & $3$ & $4$ & $5$ & $5$ & $6$ & $7$ & $7$ & $8$ & $9$ & $9$ \\
                    \hline
 \end{tabular}
 \caption{The plum-blossom wedge product table for $c=7$}
\label{tab:jiandaoji}
\end{table}

From the table for $c=7$ we see the following statements hold.
\begin{enumerate}[(1)]
\item for any $a,b$ with that $b\le 6$,  $(a,b+3)\bowtie 7=(a,b)\bowtie 7+2$;
\item for any $a,b$ with that $a\le 6$, $(a+3,b)\bowtie 7=(a,b)\bowtie 7+1$;
\item for any $a\le 9$ and any $k=1,2,3$, $(a,3k-1)\bowtie 7=(a,3k)\bowtie 7$.
\end{enumerate}

\begin{table}[!htb]
\centering
 \begin{tabular} 
 {|c||c|c|c|c|c|c|c|c|c|c|c|}
                    \hline
                    $\bowtie(c=8)$ & \rule{8pt}{0pt}$0$\rule{8pt}{0pt} &  \rule{8pt}{0pt}$1$\rule{8pt}{0pt}&  \rule{8pt}{0pt}$2$\rule{8pt}{0pt}&  \rule{8pt}{0pt}$3$\rule{8pt}{0pt}&  \rule{8pt}{0pt}$4$\rule{8pt}{0pt}&  \rule{8pt}{0pt}$5$\rule{8pt}{0pt}&  \rule{8pt}{0pt}$6$\rule{8pt}{0pt}&  \rule{8pt}{0pt}$7$\rule{8pt}{0pt}&  \rule{8pt}{0pt}$8$\rule{8pt}{0pt}&  \rule{8pt}{0pt}$9$\rule{8pt}{0pt}\\
                                   \hline\hline
                   $0$ & $0$ & $1$ & $2$ & $3$ & $3$ & $4$ & $5$ & $6$ & $7$ & $7$ \\
                   \hline
                   $1$ & $-2$ & $-1$ & $0$ & $1$ & $1$ & $2$ & $3$ & $4$ & $5$ & $5$ \\
                   \hline
                   $2$ & $-4$ & $-3$ & $-2$ & $-1$ & $-1$ & $0$ & $1$ & $2$ & $3$ & $3$ \\
                   \hline
                   $3$ & $-6$ & $-5$ & $-4$ & $-3$ & $-3$ & $-2$ & $-1$ & $0$ & $1$ & $1$ \\
                   \hline
                   $4$ & $2$ & $3$ & $4$ & $5$ & $5$ & $6$ & $7$ & $8$ & $9$ & $9$ \\
                   \hline
                   $5$ & $0$ & $1$ & $2$ & $3$ & $3$ & $4$ & $5$ & $6$ & $7$ & $7$ \\
                   \hline
                   $6$ & $-2$ & $-1$ & $0$ & $1$ & $1$ & $2$ & $3$ & $4$ & $5$ & $5$ \\
                   \hline
                   $7$ & $-4$ & $-3$ & $-2$ & $-1$ & $-1$ & $0$ & $1$ & $2$ & $3$ & $3$ \\
                   \hline
                   $8$ & $-6$ & $-5$ & $-4$ & $-3$ & $-3$ & $-2$ & $-1$ & $0$ & $1$ & $1$ \\
                   \hline
                   $9$ & $2$ & $3$ & $4$ & $5$ & $5$ & $6$ & $7$ & $8$ & $9$ & $9$ \\
                    \hline
 \end{tabular}
 \caption{The plum-blossom wedge product table for $c=8$}
\label{tab:jiandaoji}
\end{table}

From the table for $c=8$ we see that for all $a\le 9$, $(a,3)\bowtie 8=(a,4)\bowtie 8$ and  $(a,8)\bowtie 8=(a,9)\bowtie 8$.

\begin{table}[!htb]
\centering
 \begin{tabular} 
 {|c||c|c|c|c|c|c|c|c|c|c|c|}
                    \hline
                    $\bowtie(c=9)$ & \rule{8pt}{0pt}$0$\rule{8pt}{0pt} &  \rule{8pt}{0pt}$1$\rule{8pt}{0pt}&  \rule{8pt}{0pt}$2$\rule{8pt}{0pt}&  \rule{8pt}{0pt}$3$\rule{8pt}{0pt}&  \rule{8pt}{0pt}$4$\rule{8pt}{0pt}&  \rule{8pt}{0pt}$5$\rule{8pt}{0pt}&  \rule{8pt}{0pt}$6$\rule{8pt}{0pt}&  \rule{8pt}{0pt}$7$\rule{8pt}{0pt}&  \rule{8pt}{0pt}$8$\rule{8pt}{0pt}&  \rule{8pt}{0pt}$9$\rule{8pt}{0pt}\\
                                   \hline\hline
                   $0$ & $0$ & $1$ & $2$ & $3$ & $4$ & $5$ & $6$ & $6$ & $7$ & $8$ \\
                   \hline
                   $1$ & $-1$ & $0$ & $1$ & $2$ & $3$ & $4$ & $5$ & $5$ & $6$ & $7$ \\
                   \hline
                   $2$ & $-2$ & $-1$ & $0$ & $1$ & $2$ & $3$ & $4$ & $4$ & $5$ & $6$ \\
                   \hline
                   $3$ & $-3$ & $-2$ & $-1$ & $0$ & $1$ & $2$ & $3$ & $3$ & $4$ & $5$ \\
                   \hline
                   $4$ & $-4$ & $-3$ & $-2$ & $-1$ & $0$ & $1$ & $2$ & $2$ & $3$ & $4$ \\
                   \hline
                   $5$ & $-5$ & $-4$ & $-3$ & $-2$ & $-1$ & $0$ & $1$ & $1$ & $2$ & $3$ \\
                   \hline
                   $6$ & $-6$ & $-5$ & $-4$ & $-3$ & $-2$ & $-1$ & $0$ & $0$ & $1$ & $2$ \\
                   \hline
                   $7$ & $3$ & $4$ & $5$ & $6$ & $7$ & $8$ & $9$ & $9$ & $10$ & $11$ \\
                   \hline
                   $8$ & $2$ & $3$ & $4$ & $5$ & $6$ & $7$ & $8$ & $8$ & $9$ & $10$ \\
                   \hline
                   $9$ & $1$ &  $2$ & $3$ & $4$ & $5$ & $6$ & $7$ & $7$ & $8$ & $9$  \\
                    \hline
 \end{tabular}
 \caption{The plum-blossom wedge product table for $c=9$}
\label{tab:jiandaoji}
\end{table}

From the table for $c=9$ we see the following statements hold.
\begin{enumerate}[(1)]
\item for any $a\le 9$,  $(a,6)\bowtie 9=(a,7)\bowtie 9$;
\item for any $a,b$ with that $a\le 6$ and with that $b\le 6$, $(a,b)\bowtie 9=b-a$;
\item for any $a,b$ with that $a\le 6$ and with that $b\ge 7$, $(a,b)\bowtie 9=b-a-1$;
\item for any $a,b$ with that $a\ge 7$ and with that $b\ge 7$, $(a,b)\bowtie 9=9+(b-a)$.
\end{enumerate}

\begin{table}[!htb]
\centering
 \begin{tabular} 
 {|c||c|c|c|c|c|c|c|c|c|c|c|}
                    \hline
                    $\bowtie(c=1)$ & \rule{8pt}{0pt}$0$\rule{8pt}{0pt} &  \rule{8pt}{0pt}$1$\rule{8pt}{0pt}&  \rule{8pt}{0pt}$2$\rule{8pt}{0pt}&  \rule{8pt}{0pt}$3$\rule{8pt}{0pt}&  \rule{8pt}{0pt}$4$\rule{8pt}{0pt}&  \rule{8pt}{0pt}$5$\rule{8pt}{0pt}&  \rule{8pt}{0pt}$6$\rule{8pt}{0pt}&  \rule{8pt}{0pt}$7$\rule{8pt}{0pt}&  \rule{8pt}{0pt}$8$\rule{8pt}{0pt}&  \rule{8pt}{0pt}$9$\rule{8pt}{0pt}\\
                                   \hline\hline
                   $0$ & $0$ & $0$ & $0$ & $0$ & $1$ & $1$ & $1$ & $1$ & $1$ & $1$ \\
                   \hline
                   $1$ & $1$ & $1$ & $1$ & $1$ & $2$ & $2$ & $2$ & $2$ & $2$ & $2$ \\
                   \hline
                   $2$ & $2$ & $2$ & $2$ & $2$ & $3$ & $3$ & $3$ & $3$ & $3$ & $3$ \\
                   \hline
                   $3$ & $3$ & $3$ & $3$ & $3$ & $4$ & $4$ & $4$ & $4$ & $4$ & $4$ \\
                   \hline
                   $4$ & $-6$ & $-6$ & $-6$ & $-6$ & $-5$ & $-5$ & $-5$ & $-5$ & $-5$ & $-5$ \\
                   \hline
                   $5$ & $-5$ & $-5$ & $-5$ & $-5$ & $-4$ & $-4$ & $-4$ & $-4$ & $-4$ & $-4$ \\
                   \hline
                   $6$ & $-4$ & $-4$ & $-4$ & $-4$ & $-3$ & $-3$ & $-3$ & $-3$ & $-3$ & $-3$ \\
                   \hline
                   $7$ & $-3$ & $-3$ & $-3$ & $-3$ & $-2$ & $-2$ & $-2$ & $-2$ & $-2$ & $-2$ \\
                   \hline
                   $8$ & $-2$ & $-2$ & $-2$ & $-2$ & $-1$ & $-1$ & $-1$ & $-1$ & $-1$ & $-1$ \\
                   \hline
                   $9$ & $-1$ & $-1$ & $-1$ & $-1$ & $0$ & $0$ & $0$ & $0$ & $0$ & $0$ \\
                    \hline
 \end{tabular}
 \caption{The plum-blossom wedge product table for $c=1$}
\label{tab:jiandaoji}
\end{table}

From the table for $c=1$ we see that
for any $a$, $(a,0)\bowtie 2=(a,1)\bowtie 2=(a,2)\bowtie 2=(a,3)\bowtie 2$,
$(a,4)\bowtie 2=(a,5)\bowtie 2=(a,6)\bowtie 2=(a,7)\bowtie 2=(a,8)\bowtie 2=(a,9)\bowtie 2$.

\section{The plum-blossom wedge product method of multiplication of multi-digit numbers}\label{wedge-product-multiplication}

The importance of the concept of plum-blossom wedge products lies in that it can be applied to mental multiplication and mental division of multi-digit numbers so that we obtain the new method called plum-blossom wedge product methods. In this section, we introduce the  plum-blossom wedge product method of multiplication of multi-digit numbers by means of some examples.

For any multidigit multiplicand and any units multiplier, their multiplication can be mentally accomplished by a serious of concessive wedge products. For example, for computing $35649758\times 9$, one may compute in mind the following concessive wedge products:
$03\bowtie 9$, $35\bowtie 9$, $56\bowtie 9$, $64\bowtie 9$, $49\bowtie 9$, $97\bowtie 9$, $75\bowtie 9$, $58\bowtie 9$, $80\bowtie 9$,
to obtain the sequence $$(3-0,\, 5-3,\, 6-5,\, 4-6,\, 9-4-1,\, 9+(7-9),\, 3+5,\, 8-5-1,\, 2+0).$$
ie.,
$$(3,\, 2,\, 1,\, -2,\, 4,\, 7,\, 8,\, 2,\, 2),$$
which means that $$35649758\times 9=320\,847\,822.$$

To corporate the basic formula of rapid multiplication with the plum-bossom wedge products, we can carry out the multiplication of any two multidigit numbers in mind. For example,
\begin{gather*}
348\times 697\\
=(\begin{bmatrix} 03 \\6\end{bmatrix},\begin{bmatrix} 034 \\69\end{bmatrix},\begin{bmatrix} 0348 \\697\end{bmatrix},\begin{bmatrix} 3480 \\697 \end{bmatrix},
\begin{bmatrix} 480 \\97 \end{bmatrix},\begin{bmatrix} 80 \\7\end{bmatrix})_{\bowtie}\\
=(03\bowtie 6,03\bowtie 9+34\bowtie 6,03\bowtie 7+34\bowtie 9+48\bowtie 6,\\
34\bowtie 7+48\bowtie 9+80\bowtie 6,48\bowtie 7+80\bowtie 9,80\bowtie 7)\\
=(2,3+1,2+1-1,4+3-2,4+2,-4)\\
=(2,4,2,5,6,-4)\\
=242\,556\\
\end{gather*}
So we obtain that $348\times 697=242\,556.$

\section{The plum-blossom wedge product method of division of multi-digit numbers}\label{wedge-product-division}

The concept of plum-blossom wedge products  can also be applied to mental division of multi-digit numbers so that we obtain the plum-blossom wedge product method of division.
\par
Note that the division is the inverse operation of multiplication. If $A\times B=C$ then $C\div B=A$. If $A\times B+r=C$ and $0\le r< B$, then
$$C\div B=A\cdots\cdots r.$$
The plum-blossom wedge product method of division is similar to the plum-blossom product method of division. The only difference between both methods lies in changing all plum-blossom products into plum-blossom wedge products. Let us see an example.

How to calculate $242558\div 697$ by plum-blossom wedge product method? We adopt the flowing vertical algorithm, which is similar the usual ones but has some differences.

\begin{table}[!htb]
\centering
 \begin{tabular} 
 { c c c c c c c c c c}
                     &  \rule{8pt}{0pt}\rule{8pt}{0pt} &  \rule{8pt}{0pt}\rule{8pt}{0pt}&  \rule{8pt}{0pt}\rule{8pt}{0pt}&  \rule{8pt}{0pt}\rule{8pt}{0pt}
                    &  \rule{8pt}{0pt}$3$\rule{8pt}{0pt}&  \rule{8pt}{0pt}$4$\rule{8pt}{0pt}&  \rule{8pt}{0pt}$8$\rule{8pt}{0pt}
                    &  \rule{8pt}{0pt}\rule{8pt}{0pt}&  \rule{8pt}{0pt}\rule{8pt}{0pt}\\
                                   \hline
                   $6$ & $9$ & $7$ & $\bowtie$ & $2$ & $ 4$ & $2$ & $5$ & $5$  & $8$ \\ \hline
                                              &&&& $2$ & $ 1$  &&&& \\ \hline
                                                 & &&&& $3$ & $2$  &&& \\ \hline
                                                 & &&&&& $-1$  &&& \\ \hline
                                                 & &&&& $3$ & $3$  &&& \\ \hline
                                                 & &&&& $2$ & $8$  &&& \\ \hline
                                                 & &&&& & $5$ & $5$  && \\ \hline
                                                 & &&&& &  & $0$  && \\ \hline
                                                 & &&&& & $5$ & $5$  && \\ \hline
                                                  & &&&& & $5$ & $5$  && \\ \hline
                                                  & &&&& &  & $0$  & $5$ & \\ \hline
                                                 & &&&& &  & &  $6$ & \\ \hline
                                                 & &&&& &  & &  $-1$ & $8$  \\ \hline
                                                 & &&&& &  & &   & $-4$\\ \hline
                                                 & &&&& &  & &  $0$ & $2$\\
 \end{tabular}
 \caption{An example of the plum-blossom wedge product method of division}
\label{tab:jiandaoji}
\end{table}
Therefore, $242558\div 697=348\cdots\cdots 2$.

Let us explain the procedure of the above calculation. In the table, the first row is the quotient and the last row is the remainder, while the second row is the divisor and the dividend, separated by the symbol $\bowtie$.
At the  first step, we determine the first digit 3 of the quotient. Then we calculate $3\times 6+J(3\clubsuit 9)$ to obtain 21, which is in the third row of the table. $24-21=3$, which is in the fourth row. In the same row, 3 is followed by 2, which comes from dividend.
$97\bowtie 3=-1$, which is in the fifth row. $32-(-1)=33$, thus $33$ is in the sixth row.
\par
Now we can determine the second digit 4 of the quotient. Then we calculate $4\times 6+J(4\clubsuit 9)$ to obtain 28, which is at the seventh row of the table.
$33-28=5$, thus $5$ is in the eighth row. At the same row, the second 5 is the same as in the dividend.  Compute $97\bowtie 4+70\bowtie 3$ to obtain $0$ in the ninth row. $55-0=55$, thus $55$ is in the tenth row.
\par
Now we can determine the third row $8$. Then we calculate $6\times 8+J(8\clubsuit 9)$ to obtain $55$ in the eleventh row. $55-55=0$, so we have $0$ in the next row. Following this $0$ is $5$, which comes from the dividend. We calculate $97\bowtie 8+70\bowtie 4$ to obtain the next row $6$. $5-6=-1$, which is followed by $8$, the last digit of the dividend.
\par
The last digits of the divisor and of the dividend are $7$ and $8$ respectively.  So we compute $70\bowtie 8$ to obtain $-4$. Consequently,
$(-1,8)-(-4)=(0,2)$, which gives the remainder in the last row.


\section{Conclusion}\label{conclusion}

Studying multiplication and division of long integers is of great importance, not only to mental arithmetic but also to computer science.
In this paper, based on the plum-bossom product method, a novel concept--the plum-blossom wedge product is proposed.
As the same as the plum-blossom product, the plum-blossom wedge product has also some very good properties.
This new concept can be applied to multiplication and to division of long integers as well. So we obtain the new methods of multiplication and division, called the plum-blossom wedge product method. This paper presents the efficiency of this new method in arithmetic. In the next research, we will apply the  plum-blossom wedge product method to computer science, especially, to designing new multipliers and dividers suitable for long integers.


\end{sloppypar}  


\begin{thebibliography}{99}
\bibitem{Bansal}
Bansal Y,  Madhu C. A novel high-speed approach for $16\times 16$ Vedic multiplication with compressor adders[J]. Computers \& Electrical Engineering, 2016, 49:39-49.
\bibitem{Rafferty}
Ciara Rafferty, M\'{a}ire O'Neill, Neil Hanley. Evaluation of large integer multiplication methods on hardware. IEEE Trans. Comput. 66 (2017), no. 8, 1369--1382. https://ieeexplore.ieee.org/document/7869256.
\bibitem{Harvey-Hoeven}  David Harvey, Joris van der Hoeven. Integer multiplication in time $O(n \log n)$. Annals of Mathematics,
Princeton University, Department of Mathematics, In press. ffhal-02070778v2.
https://hal.archives-ouvertes.fr/hal-02070778v2.
\bibitem{Garg}
Garg A,  Joshi G. Gate diffusion input based 4-bit Vedic multiplier design[J]. IET Circuits Devices \& Systems, 2018, 12(6).
\bibitem{Gupta}
Gupta M A,  Malviya M U,  Kapse P V. A novel approach to design high speed arithmetic logic unit based on ancient Vedic multiplication technique.  2014.
\bibitem{Gurumurthy}
Gurumurthy K S,  Prahalad M S. Fast and power efficient $16\times 16$ array of array multiplier using Vedic multiplication[J]. IEEE, 2010.
\bibitem{Sadofsky} Hal Sadofsky. Fast multiplication. 2010.
http://noether.uoregon.edu/sadofsky/fast-multiplication.pdf. or
https://www.docin.com/p-1293407924.html.
\bibitem{San}
I San,  N At. On increasing the computational efficiency of long integer multiplication on FPGA[J].  Proc. 11th IEEE Int. Conf. Trust Secur. Privacy Comput. Commun., pp. 1149-1154, 2012.
\bibitem{Kavita}
Kavita U G. Performance analysis of various Vedic techniques for multiplication[J]. International Journal of Engineering Trends \& Technology, 2013, 4(3).
\bibitem{Kayal}
Kayal D,  Mostafa P,  Dandapat A, et al. Design of high performance 8 bit multiplier using Vedic multiplication algorithm with McCMOS technique[J]. Journal of Signal Processing Systems, 2014, 76(1):1-9.
\bibitem{Lee}
Lee S, Cho S M, Kim H,  et al. A combined single trace attack on global shuffling long integer multiplication and its novel countermeasure[J]. IEEE Access, vol. 8, pp. 5244-5255, 2020, doi: 10.1109/ACCESS.2019.2963317. 
\bibitem{Mathur}
Mathur M, Aarnav. Demystification of Vedic multiplication algorithm[J]. American Journal of Computational Mathematics, 2017, 07(1):94-101.
\bibitem{Paramasivam}
Paramasivam,  Sabeenian R. An efficient bit reduction binary multiplication algorithm using Vedic methods[C]// Advance Computing Conference. IEEE, 2010.
\bibitem{Pradhan}
Pradhan M,  Panda R,  Sahu S K. MAC implementation using Vedic multiplication algorithm[J]. International Journal of Computer Applications, 2011, 21(7).
\bibitem{Ramalatha}
Ramalatha M,  Dayalan K D,  Dharani P, et al. High speed energy efficient ALU design using Vedic multiplication techniques[C]// International Conference on Advances in Computational Tools for Engineering Applications. IEEE, 2009.
\bibitem{Rashno}
Rashno M,  Haghparast M,  Mosleh M. A new design of a low-power reversible Vedic multiplier[J]. International Journal of Quantum Information, 2020, 18(5):2050002.
\bibitem{Sahu}
Sahu S R,  Bhoi B K,  Pradhan M. Fast signed multiplier using Vedic Nikhilam algorithm[J]. IET Circuits Devices \& Systems, 2020, 14(8):1160-1166.
\bibitem{Shi-Wang}
Shi W,  Wang H P,  Choy C S, et al. A 0.35 V 376 Mb/s configurable long integer multiplier for subthreshold encryption[J]. IEEE Transactions on Circuits \& Systems II Express Briefs, 2018, PP:1-1.
\bibitem{Hua}
Siliang Hua, Huiguo Zhang, Jingya Zhang, Shuchang Wang. Optimization and implementation of the number theoretic transform butterfly unit for large integer multiplication[J]. Journal of Information Security and Applications, Volume 59, 2021, 102857. https://doi.org/10.1016/j.jisa.2021.102857. or \\ https://www.sciencedirect.com/science/article/pii/S2214212621000909.
\bibitem{Karatsuba}
A A Karatsuba, Y Ofman. Multiplication of multidigit numbers on automata. Soviet Physics Doklady, vol. 7, pp. 595-596, 1963,\\ http://cr.yp.to/bib/entries.html\#1963/karatsuba.
\bibitem{Shi}
Shi F S. The rapid calculation method of Shi Fengshou. Beijing: Science Press, 1989. (in Chinese)
\bibitem{Zhu}
Zhu Y W. The theory of scissor products and applications [EB/OL]. Beijing: Sciencepaper Online [2020-11-25]. \\ http://www.paper.edu.cn/releasepaper/content/202011-59254.
\bibitem{Zhu1}
Zhu Y W. The Plum-Blossom Product Method of Large Digit Multiplication and Its Application to Computer Science[J].
International Journal of Computer Applications 183(41):17-23, December 2021.\\
https://www.ijcaonline.org/archives/volume183/number41/32202-2021921805
\bibitem{Zhu2}
Zhu Y W.  A Multiplication Formula and Its Application. arXiv:2110.01820 [math.NT]. \\
https://arxiv.org/abs/2110.01820
\bibitem{Zhu3}
Zhu Y W.  On the Nine-Palace Arithmetic -- A New
Method of Mental Calculation. arXiv:2107.06647 [math.HO].\\
http://arxiv.org/abs/2107.06647
\bibitem{Zhu4}
Zhu Y W. Quick Calculation. Posts and Telecom Press. Beijing: 2023\\
\end{thebibliography}
\end{document}